\documentstyle[11pt]{article}

\input psfig

\newif\ifboxfigure      
\boxfiguretrue

\def\BoxIt#1#2{
	\vbox{\hrule
	\hbox{\vrule\kern#2\vbox{\kern#2#1\kern#2}\kern#2\vrule}
		   \hrule}}

\def\RasterBoxCap #1 #2 #3 #4 #5{

\dimen5=65pt
\divide\dimen5 by 72

\dimen0=#2\dimen5
\divide\dimen0 by 1000
\dimen1=#3\dimen5
\divide\dimen1 by 1000
\dimen2=#3\dimen5
\divide\dimen2 by 1000
\dimen3=#2\dimen5
\divide\dimen3 by 1000

\setbox4=\hbox to #4\dimen0{
 \vbox to #4\dimen1{
 \vss
 \psfig{figure=#1,height=#4\dimen2,width=#4\dimen3}
 }
 \hss
 }
 \vbox {
 \ifboxfigure\BoxIt{\box4}{0pt}
 \else\box4
 \fi
 \vskip 2pt
 \hbox to #4\dimen0{\hfil #5 \hfil }
 \vbox to 2pt{\relax}
 }
}

\begin{document}
\large
\centerline {\Large \bf Dynamics of certain non-conformal semigroups}

\vskip5pt
\centerline{Yunping Jiang }
\centerline{Institute for Mathematical Sciences} 
\centerline{SUNY at Stony Brook, Stony Brook, NY 11794}
\centerline{January, 1991, revised December, 1991}

\vskip40pt
\centerline{ {\Large \bf Abstract}}

\vskip5pt
A semigroup generated by two dimensional $C^{1+\alpha}$ 
contracting maps is considered. We call a such semigroup regular
if the maximum $K$ of the conformal dilatations of generators, the
maximum $l$ of the norms of the derivatives of generators and the
smoothness $\alpha$ of the generators satisfy a compatibility condition
$K< 1/l^{\alpha}$.  We prove that the shape of the image of the core of
a ball under any element
of a regular semigroup is good (bounded geometric distortion like the Koebe
$1/4$-lemma \cite{a}). And we use it to show a lower and a upper bounds of
the Hausdorff dimension of 
the limit set of a regular semigroup. We also consider a semigroup generated
by higher dimensional maps.

\vskip30pt
\centerline{{\Large \bf Contents}}

\vskip5pt
\noindent \S 0 Introduction.

\noindent \S 1 Statements of main results.

\noindent \S 2 Proof of Theorem 1.

\noindent \S 3 Proof of Theorem 2.

\noindent \S 4 Higher dimensional regular semigroups and some remarks.

\pagebreak
\noindent { \bf \S 0 Introduction.}

It is a well-known result \cite{r,s1} 
that the Hausdorff dimension of the
Julia set of a complex quadratic polynomial $p(z)= z^{2} +c$ is greater than 
one for a complex number $c$ with small $|c|\neq 0$ 
(see \cite{b1} for a similar result
in quasifuchsian groups). Now consider a  
non-conformal complex map $f(z)=z^{2}+b\overline{z}+c$ where $b$ and 
$c$ are complex
parameters (or 
$f(z)= z^{n}|z|^{(\gamma-n)} + c$ where 
$\gamma>0$ is a real parameter, $c$ is a complex parameter 
and $n>0$ is a fixed integer). Let $\lambda =(b, c)$ (or $\lambda =(\gamma
-n, c)$ and $|\lambda |= |b|+|c|$ (or $|\lambda |= |\gamma -n|+|c|$).
The map $f_{0}(z)=z^{2}$ (or $f_{0}(z)=z^{n})$ is analytic and expanding on a
neighborhood $U$ of $S^{1}=\{ z \in {\bf C}; |z|=1\}$ which is the
maximal invariant set of $f_{0}$ in $U$. By the structural stability theorem
\cite{sh}, for $|\lambda|$ small, there is a set $J_{\lambda}$ such that 
it is the maximal invariant set of $f$ and $f|J_{\lambda}$ is
conjugate to $f_{0}|S^{1}$, that is, there is a homeomorphism $h$ from a
neighborhood of 
$S^{1}$ onto a neighborhood of $J_{\lambda}$ such that 
$f\circ h = h\circ f_{0}$. Thus the set
$J_{\lambda}$ is a Jordan curve. It is easy to see that $J_{\lambda}$ is the
boundary of the basin $B_{\infty}=\{ z \in {\bf C}; |f^{\circ k}(z)| \mapsto \infty$ as
$k\mapsto +\infty  \}$ for $|\lambda|$ small (see Fig.~1 and Fig. 2). 
We may call $J_{\lambda}$ the Julia set of
$f$ (ref. \cite{m}).

\begin{figure}
\centerline{\hbox{
\RasterBoxCap {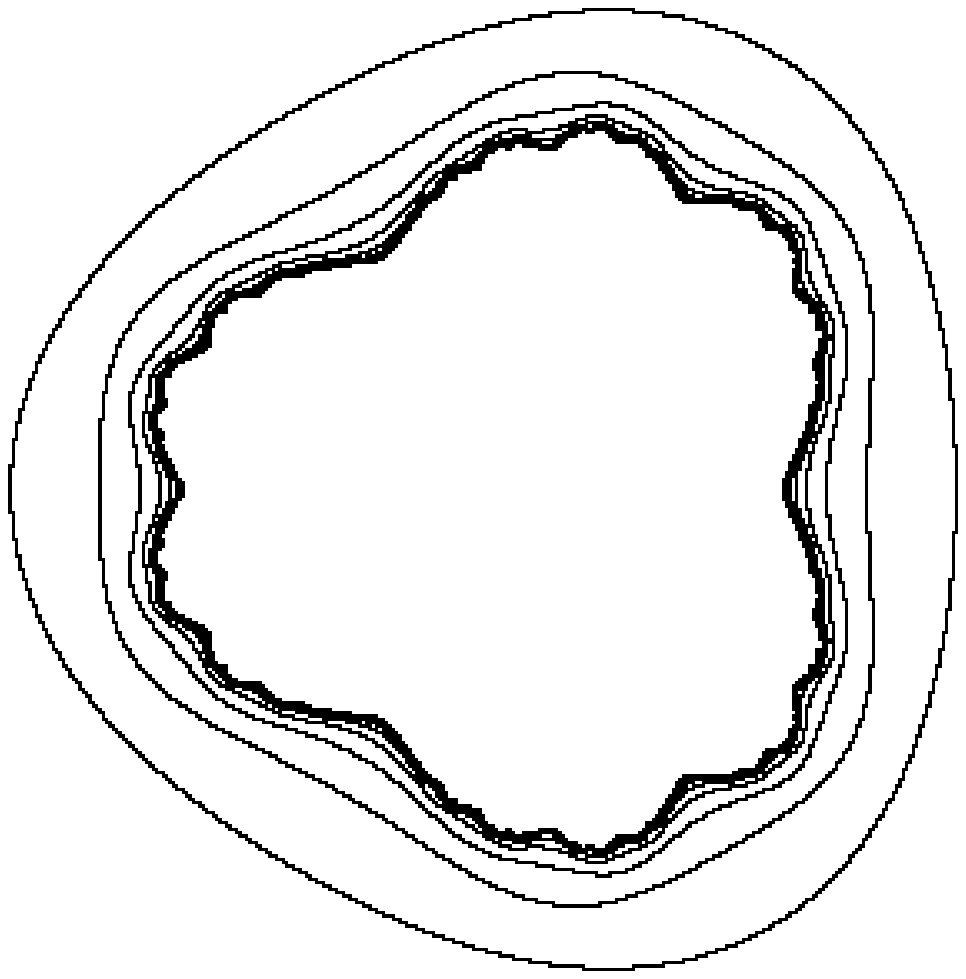} {400} {400} {300} 
{$\lambda=(0.2, 0)$}\hss
\RasterBoxCap {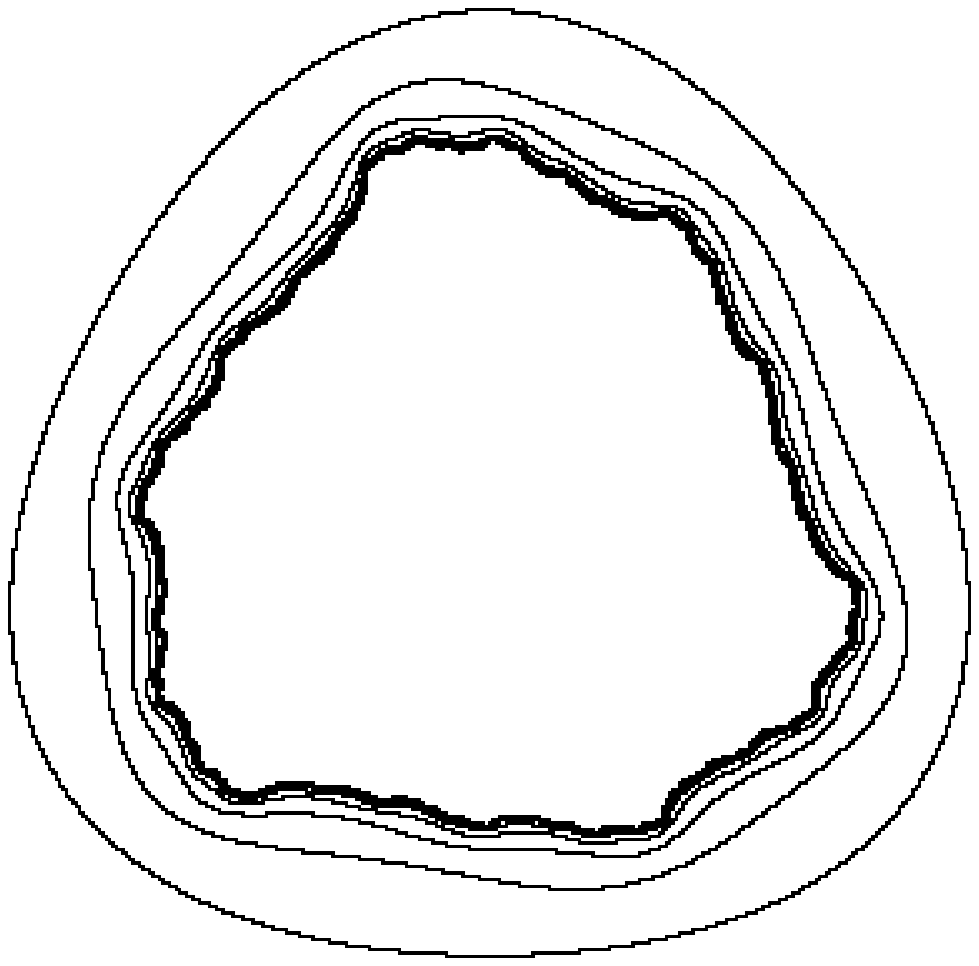} {400} {400} {300} 
{$\lambda=(0.2i, 0)$}\hss
\RasterBoxCap {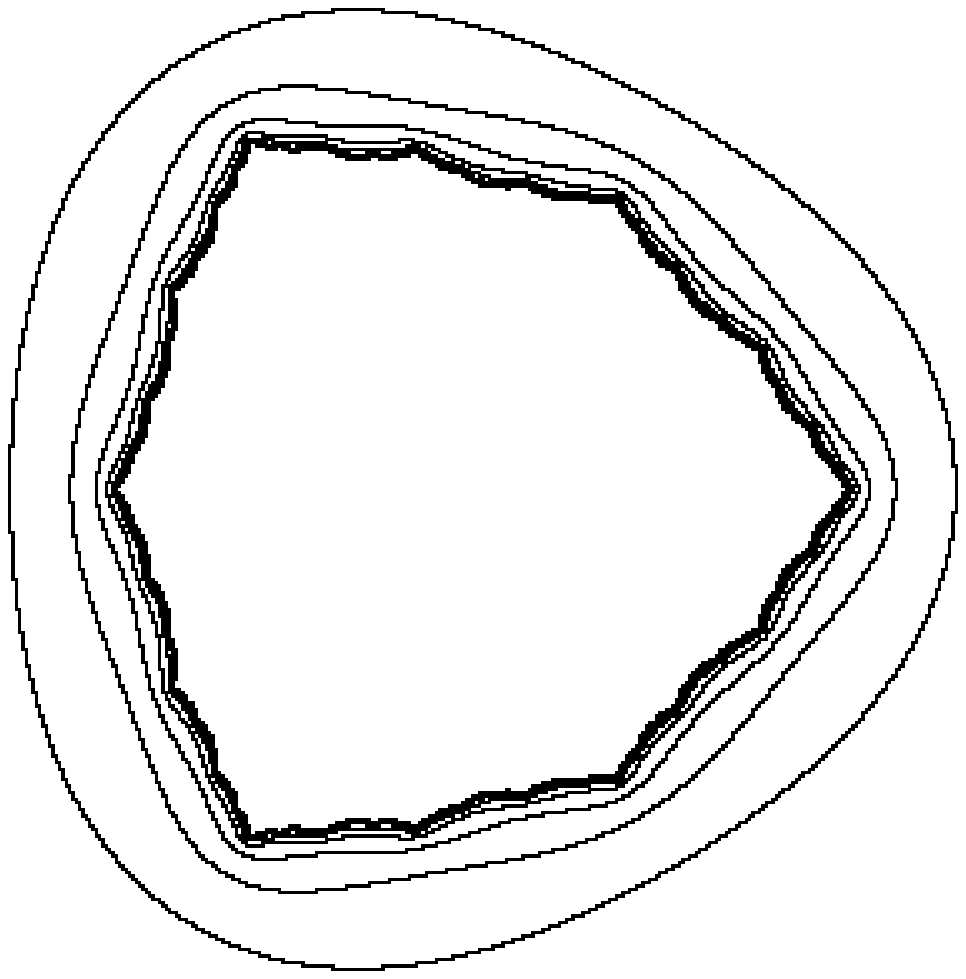} {400} {400} {300} 
{$\lambda=(-0.2, 0)$}\hss 
\RasterBoxCap {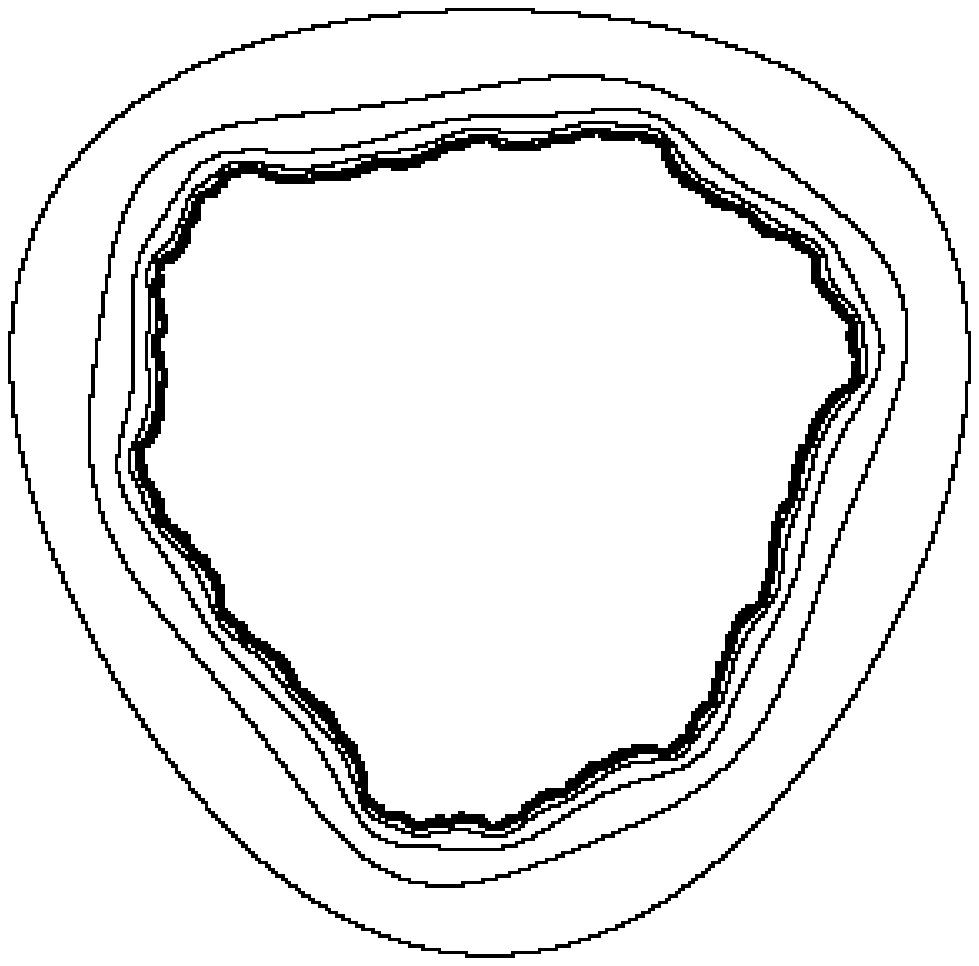} {400} {400} {300} 
{$\lambda=(-0.2i, 0)$}
}}

\centerline{\hbox{
\RasterBoxCap {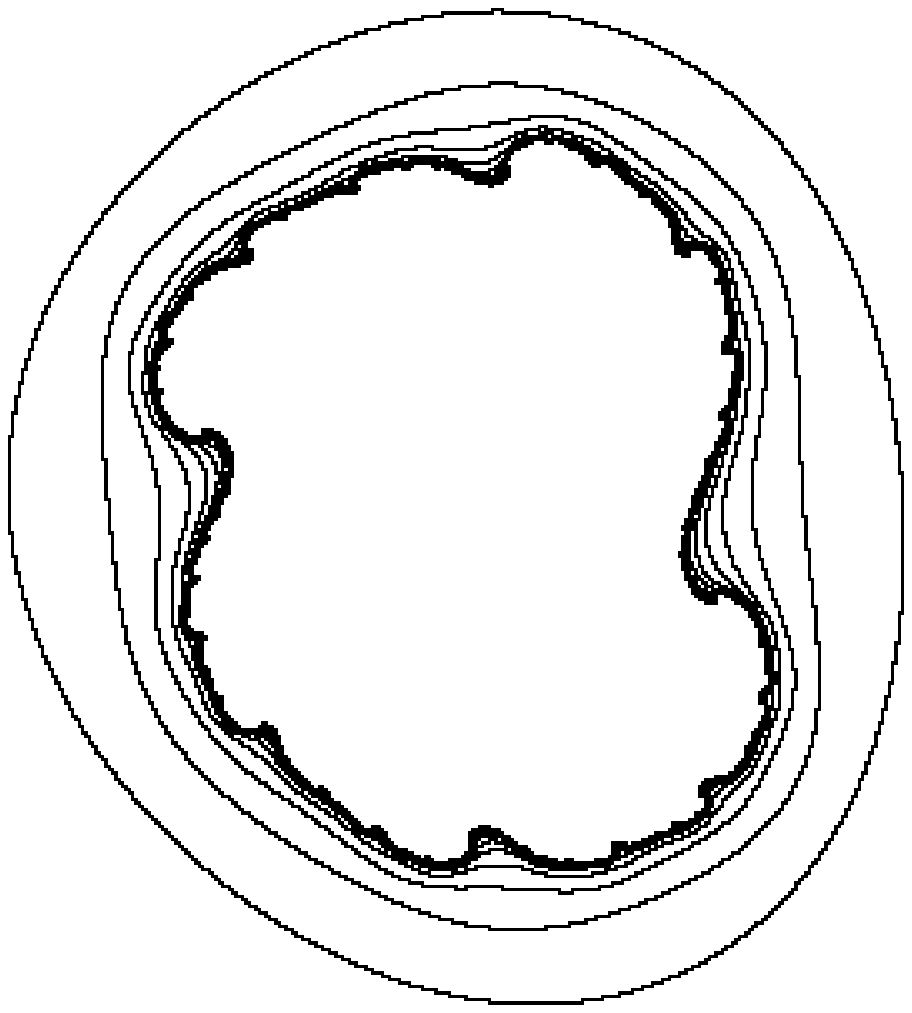} {400} {400} {300} 
{$\lambda=(0.1,0.2+0.1i)$}\hss
\RasterBoxCap {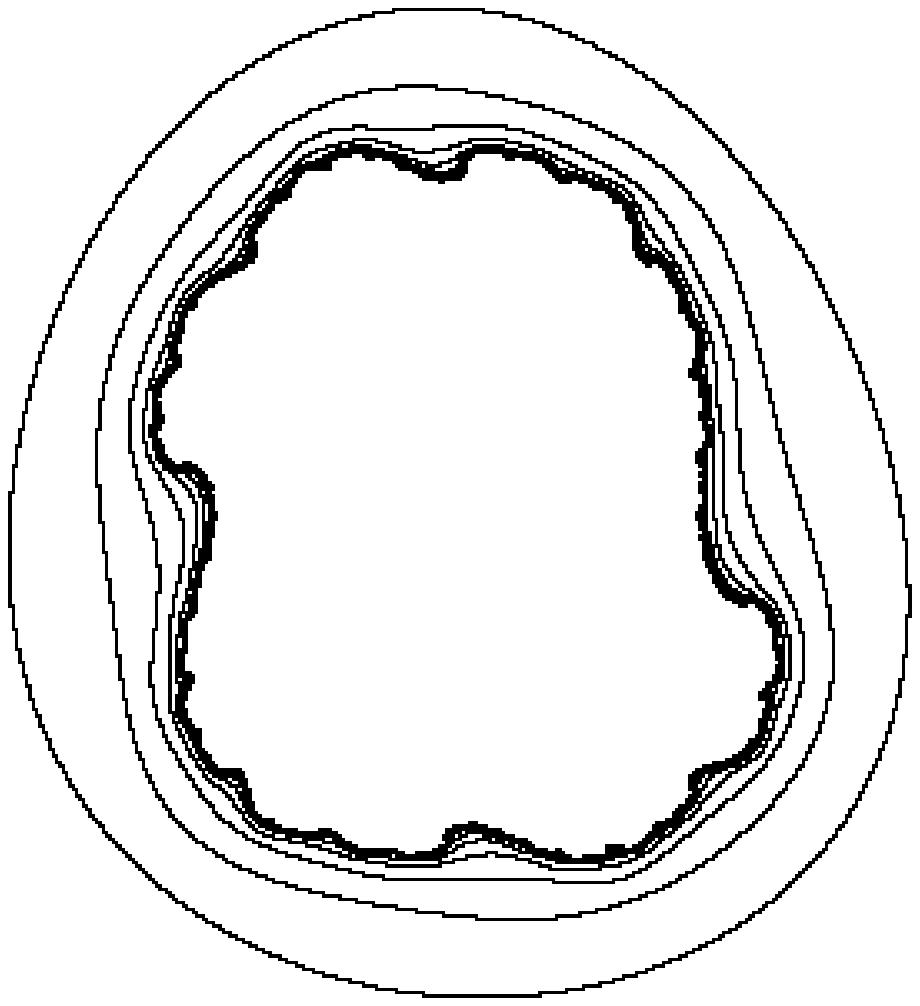} {400} {400} {300} 
{$\lambda=(0.1i,0.2+0.1i)$}\hss
\RasterBoxCap {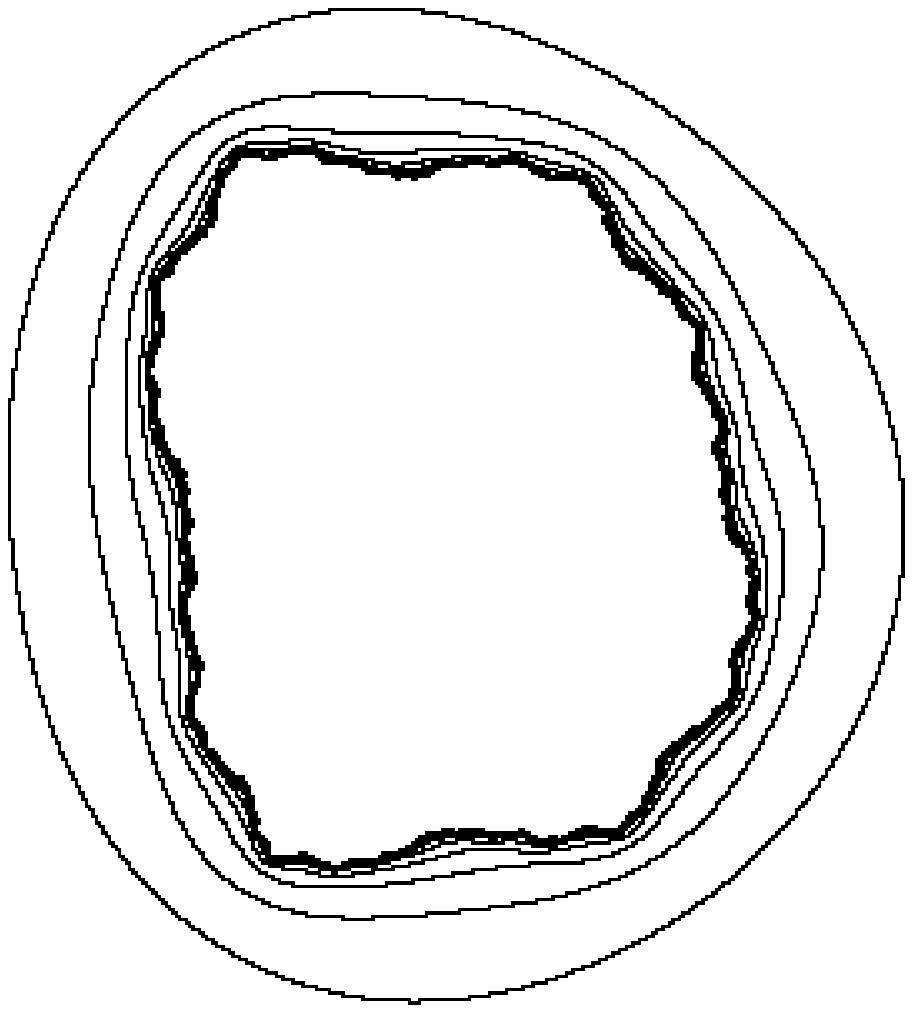} {400} {400} {300} 
{$\lambda=(-0.1,0.2+0.1i)$}\hss 
\RasterBoxCap {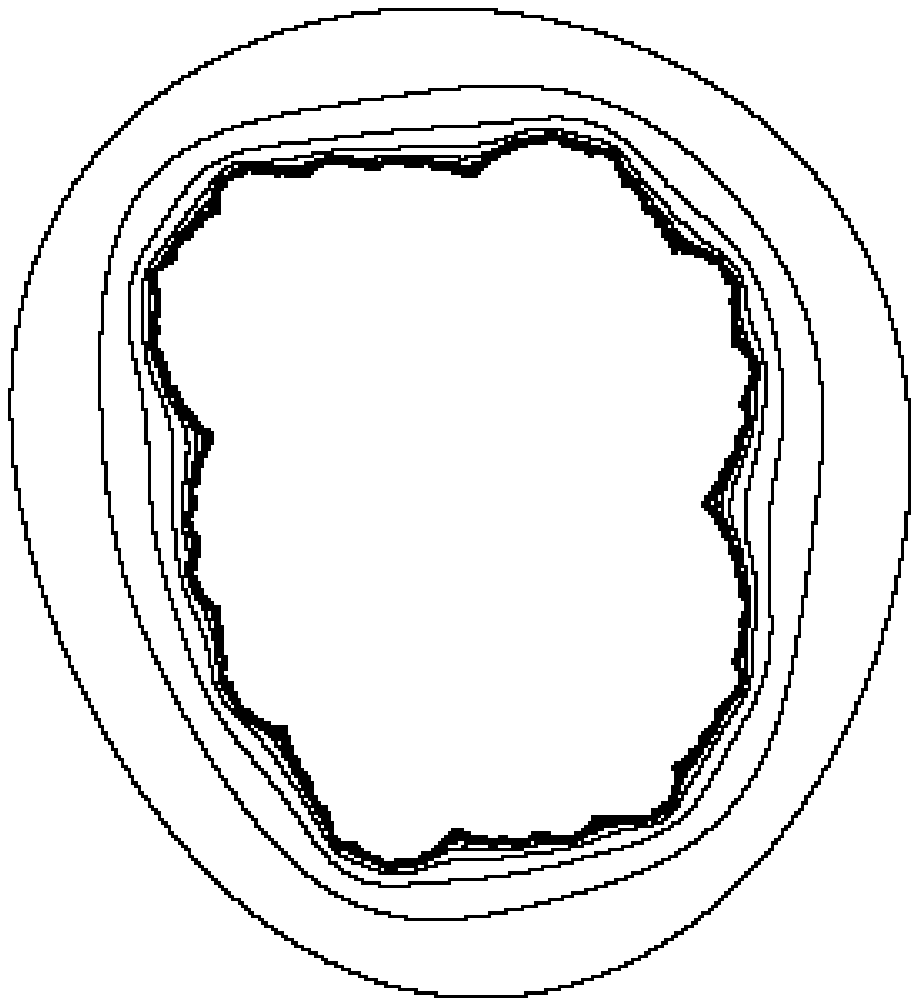} {400} {400} {300} 
{$\lambda=(-0.1i,0.2+0.1i)$}
}}

\centerline{{\bf Fig. 1:} Preimages of a circle with large radius under} 
\centerline{iterates of $f(z)=z^{2}+b\overline{z} +c$ and $\lambda =(b,c)$.}

\vskip20pt
\centerline{\hbox{
\RasterBoxCap {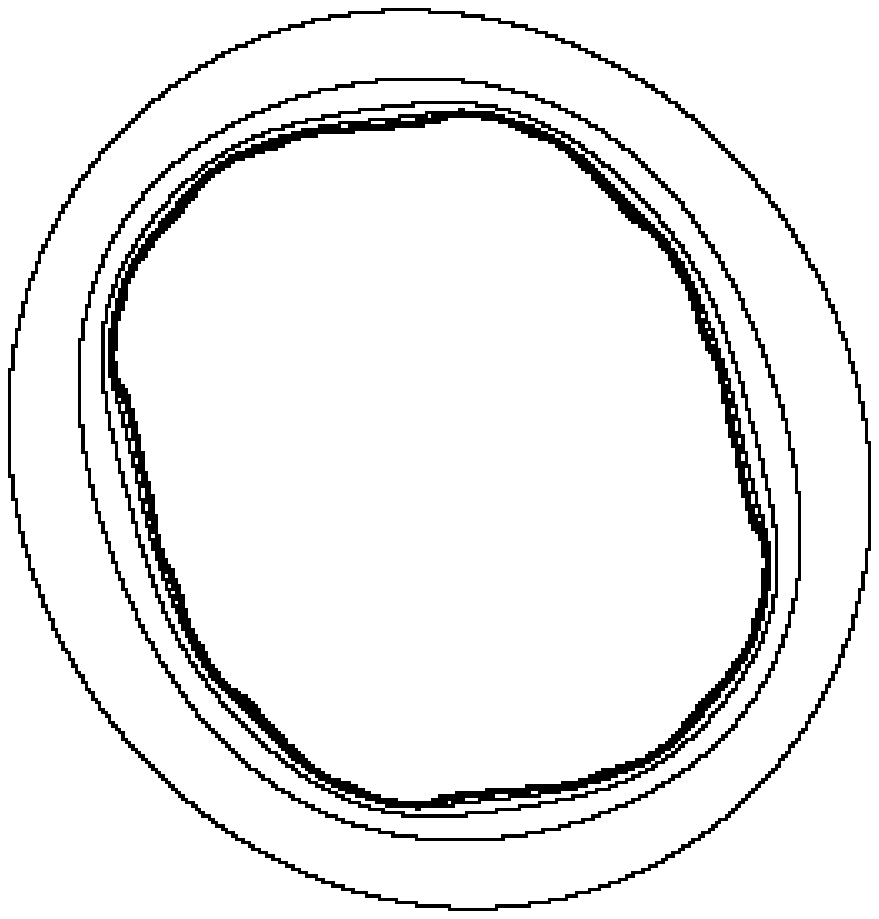} {400} {400} {300} 
{$\lambda=(0.5, 0.1+0.2i)$}\hss
\RasterBoxCap {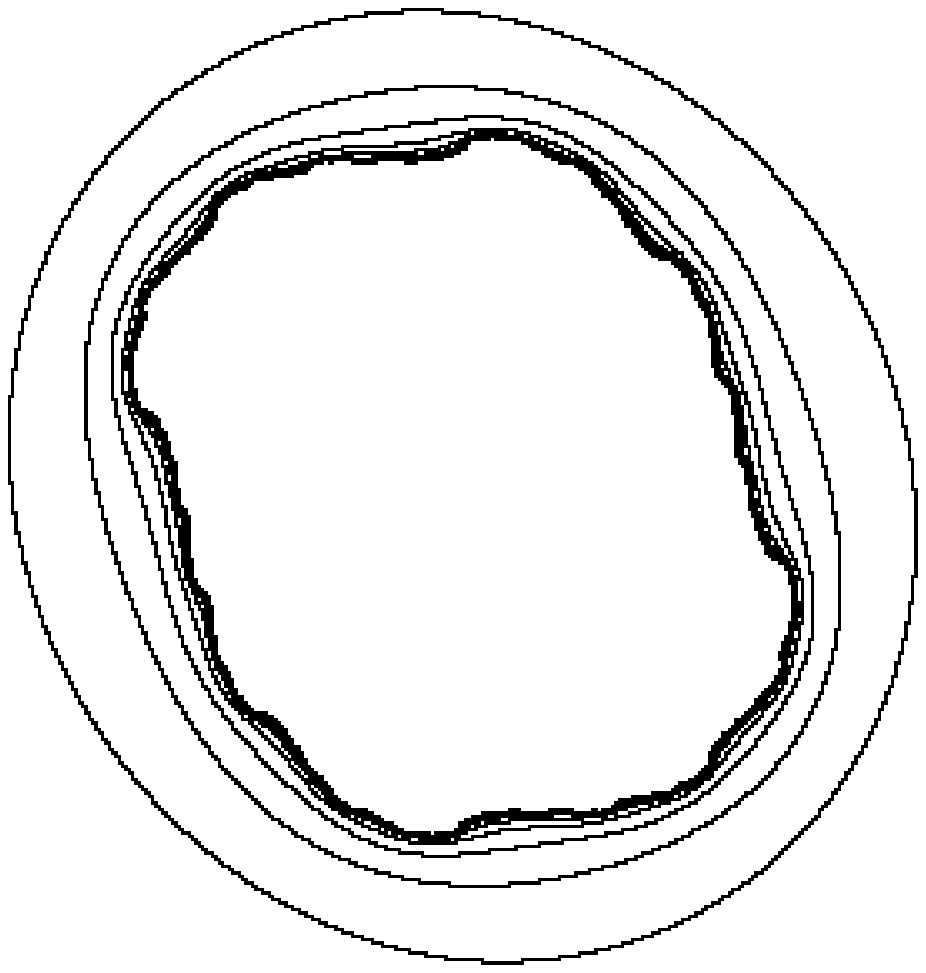} {400} {400} {300} 
{$\lambda=(0.1, 0.1+0.2i)$}\hss
\RasterBoxCap {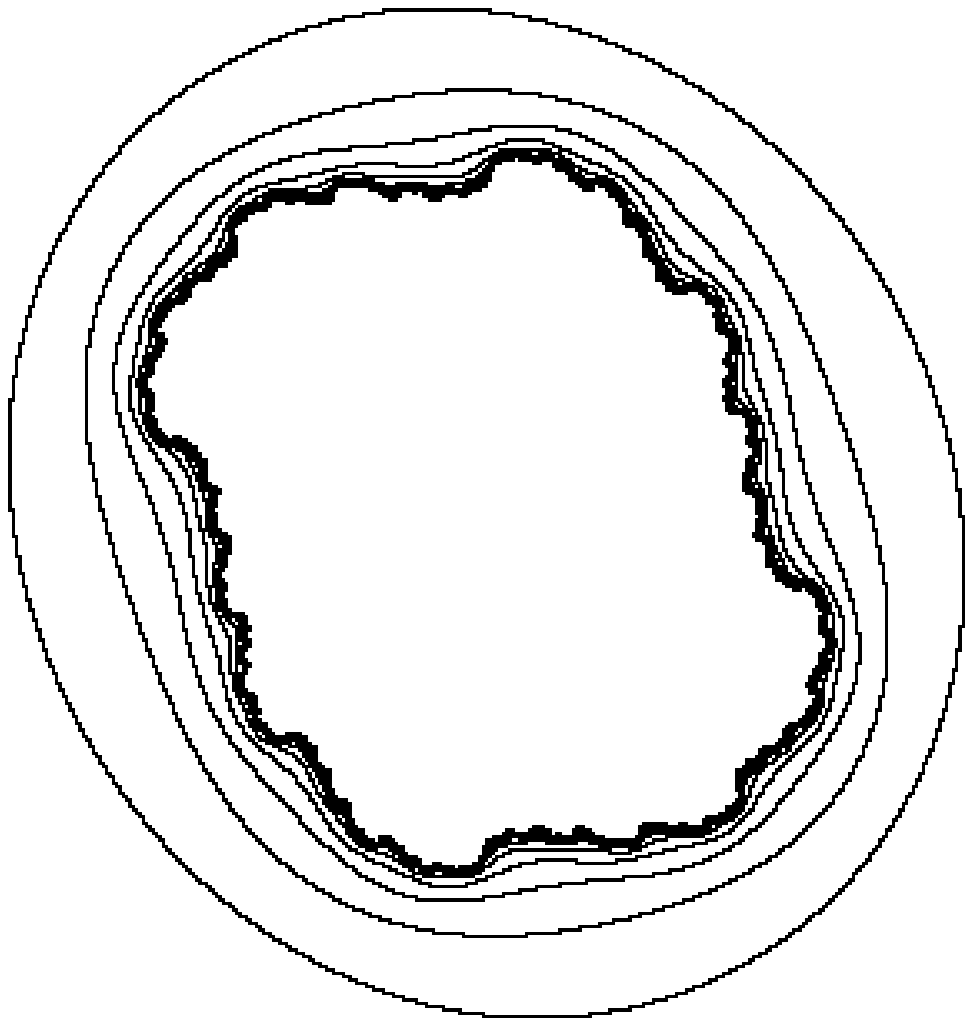} {400} {400} {300} 
{$\lambda=(-0.2, 0.1+0.2i)$}\hss 
\RasterBoxCap {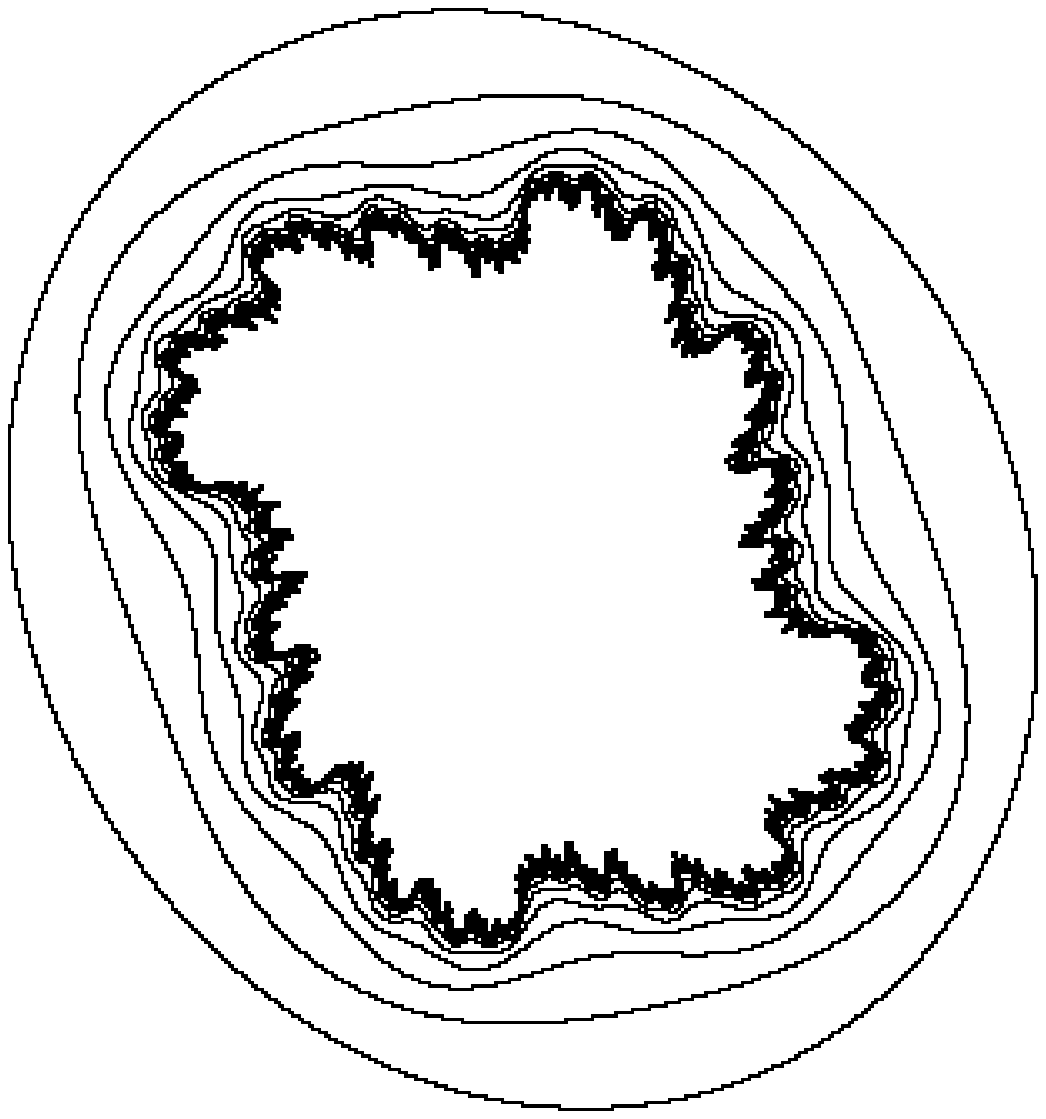} {400} {400} {300} 
{$\lambda=(-0.5, 0.1+0.2i)$}
}} 

\centerline{{\bf Fig. 2:} Preimages of a circle with large radius under}
\centerline{iterates of $f(z)=z^{2}|z|^{\gamma-2} +c$ and $\lambda =(\gamma-2, c)$.}
\end{figure}

\vskip5pt
{\bf Question 1.} {\em Is the Hausdorff dimension of the Julia set
$J_{\lambda}$ 
of $f(z)$ greater than $1$ for some small $|b| \neq 0$ and small $|c|$ (or
small $|c|\neq 0$ and 
small $|\gamma-n| \neq 0$) ? }

We will prove some
general results (Theorem 1
and Theorem 2) in \S 1, \S 2 and \S 3, which can be used to give the answer
(Corollary 3) to this question. We note that the general results 
themselves are interesting and have other applications \cite{j2}.

\vskip10pt
{\bf Acknowledgement.} I would like to thank Professor Dennis Sullivan for
very useful discussions and remarks. The conjecture in Remark 3 is formulated
when the author visited the Mathematics Institute at University of
Warwick. I would like to thank Professor David Rand for 
useful conversations. 

\vskip10pt
\noindent { \bf \S 1 Statements of main results.}

Suppose $V$ and $U$ are two bounded and open sets of the complex plane 
${\bf C}$ with $\overline{V}\subset U$ and $f$ is
a $C^{1}$-map from $U$ into ${\bf C}$. 
The restriction $f|\overline{V}$  
is said to be $C^{1+\alpha}$ for some
$0< \alpha \leq 1$ if  
\[ f(w)= f(z) + \Big (D(f)(z) \Big)(w-z) + R(w,z)\]
satisfies $| R(w,z) | \leq L_{0}|w-z|^{1+\alpha }$ for
$z\in \overline{V}$ and $w\in U$ where $L_{0}>0$ is a constant and 
$D(f)(z)$ is the derivative of $f$ at $z$.
For a $C^{1+\alpha}$ diffeomorphism $f$ from $\overline{V}$ onto $\overline{W}$, 
we use $g$ to denote
its inverse. The map $g$ is said to be contracting if 
there is a constant $0< \lambda < 1$ such that
$\left| \Big( D(g)(z)\Big) (v) \right| \leq \lambda |v|$ for all $z$ in
$\overline{W}$ and all $v$ in ${\bf C}$.
Suppose $V_{i}$ and $U_{i}$, $i=0, 1, \ldots, n-1$, are pairs of bounded
open sets of ${\bf C}$ with $\overline{V}_{i}\subset U_{i}$ and $f_{i}$
are maps from $U_{i}$ into ${\bf C}$ such that the restriction $f_{i}|\overline{V_{i}}$
from $\overline{V_{i}}$ onto $\overline{W_{i}}$ are $C^{1+\alpha}$ diffeomorphisms 
for some $0< \alpha \leq 1$ and the inverses $g_{i}$ of
$f_{i}|\overline{V_{i}}$ are
contracting. To simplify the notations, we assume 
that $W=W_{i}$ for
all $i$ and $\cup_{i=0}^{n-1}V_{i}
\subset W$.
We will use ${\cal G}=\langle g_{0}, g_{1}$, $\ldots$, $g_{n-1} \rangle$ to
denote 
the semigroup
generated by 
all $g_{i}$ and use $\Lambda
=\cap_{g\in {\cal G}}\Big( g(\overline{W})\Big)$ to denote the limit set of ${\cal G}$, which is 
compact, completely invariant (the existence of $\Lambda$ can be proven by
using Hausdorff distance on subsets).

Suppose $z=x+yi$ is a point in ${\bf C}$ and 
$\overline{z}=x-yi$ is the conjugate of $z$. By
the complex analysis \cite{a}, we know that for $z\in \overline{W}$ and
$w\in {\bf C}$ with $|w|=1$,
\[ \Big| |(g_{i})_{z}|- |(g_{i})_{\overline{z}}|\Big| \leq \left| \Big(
D(g_{i})(z)\Big)(w)\right| \leq |(g_{i})_{z}|+ |(g_{i})_{\overline{z}}|.\]
Let 
\[ l_{i}(z)= |(g_{i})_{z}|+ |(g_{i})_{\overline{z}}|, \hskip10pt
s_{i}(z)=\Big| |(g_{i})_{z}|- |(g_{i})_{\overline{z}}|\Big| \]
and $K_{i}(z)=l_{i}(z)/s_{i}(z)$, the conformal
dilatation of $g_{i}$ at $z$. Let $l=\max \{ l_{i}(z) \}<1$, 
$s=\min \{ s_{i}(z) \}>0$ and $K=\max \{ K_{i}(z) \}<+\infty$ where $\max$
and $\min$ are over
all $z$ in $\overline{W}$ and all $0\leq i <n$. 

{\bf Definition 1.} {\em We say ${\cal G}$ is 
regular if $K< 1/l^{\alpha}$.}

Denote by $B(z, r)$ the closed disk of
radius $r$ centered at $z$.
One of the main results, which generalizes the Koebe $1/4$-lemma \cite{bi} 
in some sense, is the following:

{\bf Theorem 1} (geometric distortion). {\em Suppose
${\cal G}=\langle g_{0}$, $g_{1}$, $\ldots$, $g_{n-1}\rangle $ is regular. 
There are two functions
$\delta=\delta (\varepsilon)>0$ and $C=C(\varepsilon)\geq 1$ with  
$\delta (\varepsilon) \mapsto 0$ and $C(\varepsilon)\mapsto 1$ as $\varepsilon
\mapsto 0+$ such that 
\[ g\Big(
B(z, r)\Big)\supset g(z)+C^{-1}\cdot \Big( D(g)(z)\Big) \Big( B(0, r )\Big)
\hskip7pt and \]
\[ g\Big(
B(z, r)\Big) \subset
g(z)+C\cdot \Big( D(g)(z)\Big)\Big( B(0, r )\Big) \]
for any $0< r \leq \delta (\varepsilon)$, any $g\in
{\cal G}$ and any $z\in \overline{W}$ (see Fig. 3).}

\vskip7pt
Let $\angle \Big( g(w)-g(z), \Big( D(g)(z)\Big) (w-z)\Big)$ be the
smallest angle between the vectors $g(w)-g(z)$ and $\Big( D(g)(z)\Big)
(w-z)$.

\vskip5pt
{\bf Corollary 1} (angle distortion). {\em 
Moreover, there is a function $D(\varepsilon)>0$ with $D(\varepsilon) \mapsto0$ as $\varepsilon\mapsto 0+$ such that 
\[ \left| \log \Big( \angle \Big( g(w)-g(z), \Big( D(g)(z)\Big)
(w-z)\Big) \right| \leq D(\varepsilon) \]
for $0< r \leq \delta (\varepsilon)$, $g\in
{\cal G}$, $z\in \overline{W}$ and $w\in B(z,r)$.}

\vskip-10pt
\centerline{\psfig{figure=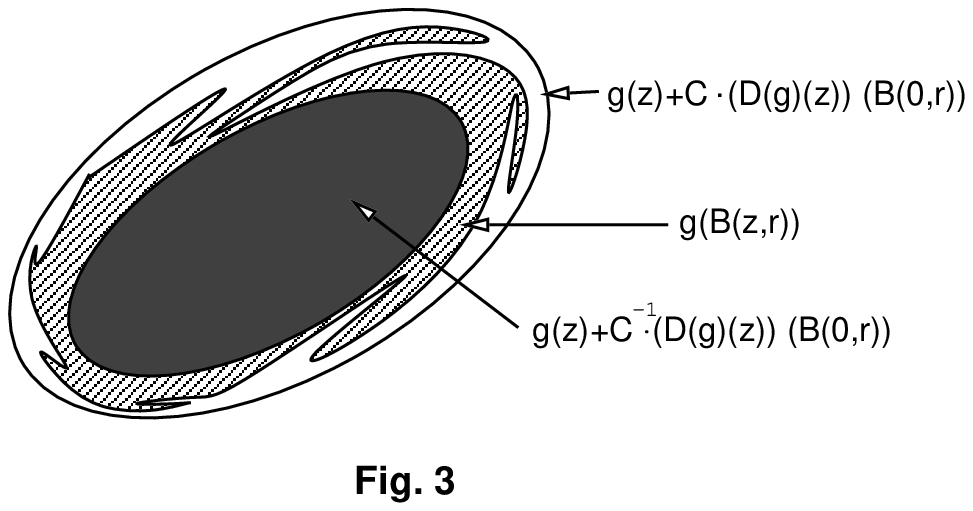}}

\vskip5pt
A regular semigroup ${\cal G}=\langle g_{0},$ $\ldots$, $g_{n-1}\rangle $ is said to be
Markov for a real number $\delta_{0}>0$ if there are simple
connected, pairwise disjoint open sets
$\Omega_{0}$, $\Omega_{1}$, $\ldots$, $\Omega_{q-1}$ such that 
\begin{itemize}
\item[(a)] $\max_{0\leq l\leq q-1}diam(\Omega_{l})\leq \delta_{0}$,
\item[(b)] $\cup_{l=0}^{q-1} \overline{\Omega_{l}} \supset \Lambda $, and 
\item[(c)] $f_{i}(\overline{\Omega_{l}\cap \Lambda}) = 
\Big( \cup_{t=1}^{k_{l}}\overline{\Omega_{i_{t}}} \Big)
\cap \Lambda$ for every $0\leq l< q$ and $\Omega_{l}\subset V_{i}$ where
$f_{i} = g_{i}^{-1}$.
\end{itemize}

\noindent Without loss of generality, we may assume  
$q=n$ and $g_{i}=(f_{i}|\Omega_{i})^{-1}$ if ${\cal G}$ is Markov. 

Suppose ${\cal G}=\langle g_{0},$ $\ldots$, $g_{n-1} \rangle$ is a regular and Markov
semigroup. Let $A=(a_{ij})$ be the $n\times n$ matrix of $0$ and $1$
such that $a_{ij}=1$ if $f_{i}(\Omega_{i}\cap \Lambda )\supset \Omega_{j}\cap
\Lambda$ and $a_{ij}=0$ 
otherwise. A sequence $w_{p}=i_{0}i_{1}\cdots i_{p-1}$ of symbols $\{ 0$, $1$,
$\ldots$, $n-1 \}$ is said to be admissible if
$a_{i_{l}i_{l+1}}=1$  for $l=0$, $1$, $\ldots$, $p-1$ ($p$ may be $\infty$).
Let $\Sigma_{p}$ be the space of all admissible sequences $w_{p}$ of
length $p$, 
$\sigma(i_{0}i_{1}\cdots )=i_{1}\cdots$ be the shift map on
$\Sigma_{\infty}$ and $\pi
(i_{0}i_{1}\cdots) =\cap_{k=0}^{\infty} g_{i_{k}}(\overline{W})$ be the
projection from $\Sigma_{\infty}$ to $\Lambda$ \cite{b,r} (note that $\pi$
is the semi-conjugacy). 
We call the functions 
\[ \phi_{up} (w) =\log \Big( l_{i}\circ
\pi (w)\Big) \hskip7pt
{\rm and}  \hskip7pt \phi_{lo} (w) =\log \Big( s_{i}\circ \pi (w)\Big), \]
for $w=ii_{1}\cdots \in \Sigma_{\infty}$,
the upper and lower potential functions of ${\cal G}$.
They are H\"older \cite{b}.

Let $P$ be the pressure function (see, for example, \cite{b,r}) defined on
$C^{H}$, the space of H\"older continuous functions on $\Sigma_{\infty}$.  
Then \cite{b}
\[ P(\phi) = \lim_{p\mapsto \infty} \frac{1}{p}\log \Big( \sum_{w\in 
fix(\sigma^{\circ p})}\exp
\Big( \sum_{k=0}^{p-1}\phi \Big( \sigma^{\circ k}(w)\Big) \Big) \Big). \] 
For $\phi = \phi_{up}$ or $\phi_{lo}$, 
$P(t\phi)$ is continuous, strictly monotone and convex 
function on the real line and tends to 
$-\infty$ and $+\infty$ as $t$ goes to
$+\infty$ and $-\infty$.
There is a unique $t_{up}>0$ ($t_{lo}>0$) such that $P(t_{up}\phi_{up})=0$ 
($P(t_{lo}\phi_{lo})=0$) \cite{b1,r}. 

{\bf Theorem 2.} {\em Suppose ${\cal G}=\langle g_{0},$ $\ldots$,
$g_{n-1}\rangle $ is a regular and Markov
semigroup and $HD(\Lambda)$ is the Hausdorff dimension of
the limit set 
$\Lambda $ of ${\cal G}$. Then $t_{lo}\leq HD(\Lambda) \leq t_{up}.$}   

Suppose ${\cal G}_{\lambda}= \langle g_{0,\lambda}$, $\ldots$,
$g_{n-1,\lambda}\rangle $ 
is a family of regular 
and Markov semigroups such that every $g_{i,\lambda}(z)$ is $C^{1}$ on both 
variables $\lambda$ and
$z$. Let $HD(\lambda)$ be the Hausdorff dimension of 
the limit set $\Lambda_{\lambda}$ of ${\cal G}_{\lambda}$.

{\bf Corollary 2.} {\em If all $g_{i,\lambda_{0}}$ are conformal
($K_{\lambda_{0}}=1$), then $HD(\lambda)$ is continuous at $\lambda_{0}$.} 

{\bf Corollary 3.} {\em Suppose 
$f(z)=z^{2}+b\overline{z}+c$ $($or f(z)=
$z^{n}|z|^{(\gamma-n)} +c )$ and $\lambda =(b,c)$ $($or $\lambda =(\gamma -n,
c))$. 
For each $c$ with small $|c|\neq 0$, there is a $\tau (c)>0$ such
that for every $|b|\leq \tau (c)$  $($ or $ |\gamma-n|\leq \tau (c) )$, 
the Hausdorff
dimension $HD(\lambda)$ of the Julia set
$J_{\lambda}$ of $f$ is bigger than one (see Fig. 4 in \S 4).}   

{\bf Remark 1.} Biefeleld, Sutherland, Tangerman and Veerman \cite{bstv}   
showed recently
that for $f(z)=z^{2}|z|^{(\gamma-2)} +c$ and 
a small $\gamma-2>0$,
there is an $\eta (\gamma) >0$ such that the Julia set $J_{\lambda}$ of $f(z)$ for $|c|< \eta
(\gamma)$ is a smooth circle (see Fig. 4 in \S 4).

\vskip10pt
\noindent {\bf \S 2 Proof of Theorem 1.}

By the compactness of $\overline{W}$, there is a 
function $\delta=\delta (\varepsilon)>0$ with $\delta (\varepsilon
)\mapsto 0$ as $\varepsilon \mapsto 0+$ such that every $g_{i}$ is defined on 
$B(z, \delta)$ for $z$ in $\overline{W}$ and $g_{i}(w) =
g_{i}(z) +\Big( D(g_{i})(z)\Big)(w-z) +R_{i}(w,z)$
satisfies that 
\[ |R_{i}(w,z)|\leq \Big( \varepsilon/2 \Big) \cdot \Big( \inf_{w\in
\overline{W}}||D(g_{i})(z)||\Big) \cdot |w-z|\] 
for $z$ and $w$ in
$\overline{W}$ with
$|w-z|\leq \delta$ and $0\leq i<n$.
This implies that for $z$ in $\overline{W}$ and $0< r\leq \delta$,
\[ g_{i}\Big( B(z, r)\Big) \supset 
g_{i}(z)+(1+\varepsilon)^{-1} \cdot \Big( D(g_{i})(z)\Big) \Big( B(0, r
)\Big) \hskip7pt and \] 
\[ g_{i}\Big( B(z, r)\Big)
\subset
g_{i}(z)+(1+\varepsilon) \cdot \Big( D(g_{i})(z)\Big) \Big( B(0, r )\Big)
\hskip30pt (*).\]

Suppose $L_{0}>0$ and $0< \beta <\alpha$ are constants 
such that $|R_{i}(w,z)|\leq
L_{0}|w-z|^{1+\alpha}$ and $K_{i}(z) \leq \Big( 1/l_{i}(z) \Big)^{\beta}$ 
for $0\leq i <n$, $z$ and $w$ in $\overline{W}$. Let 
$ \kappa_{m}=\sum_{i=0}^{m}l^{(\alpha -\beta)i}$. We take $\delta= \delta
(\varepsilon)\leq 1$ so small that 
\[ \Theta_{\varepsilon} = \Big( L_{0}/s\Big) \big( 1+\varepsilon
+\kappa_{\infty}\Big)^{1+\alpha} \delta^{(\alpha -\beta)} \leq 1 \]
and then take
\[ C_{m}(\varepsilon)= 1+\varepsilon + \delta^{\beta} \cdot \kappa_{m}.\] 
It is clear that $C_{m}(\varepsilon)
\mapsto 1$ as $\varepsilon \mapsto 0+$.
 
{\bf Claim.} For $g=g_{i_{0}}\circ g_{i_{1}}\circ \cdots \circ g_{i_{m}}$ in
${\cal G}$, \[ g\Big(
B(z, r)\Big) \supset g(z)+C_{m}^{-1}\cdot \Big( D(g)(z)\Big) \Big( B(0,
r ) \Big) \hskip7pt and \]
\[g\Big(
B(z, r)\Big)
 \subset
g(z)+C_{m}\cdot \Big( D(g)(z)\Big) \Big( B(0, r) \Big).\]

{\bf Proof of claim.} For $m=0$, it is the formulae in $(*)$. 
Suppose the claim holds for $m=0$,
$1$, $\ldots$, $M-1$ ($M\geq 1$). 
Then for $g=g_{i_{0}}\circ g_{i_{1}}\circ \cdots \circ
g_{i_{M}}=g_{i_{0}}\circ G$, 
\[ 
g\Big( B(z, r)\Big)\supset g_{i_{0}}\Big( G(z)+C_{M-1}^{-1}\cdot \Big(
D(G)(z)\Big) \Big( B(0, r )
\Big) \Big) \hskip7pt and \]
\[ g\Big( B(z, r)\Big) \subset
g_{i_{0}}\Big( G(z)+C_{M-1} \cdot \Big( D(G)(z)\Big) \Big( B(0, r )
\Big) \Big).\]
For any $w$ in $B(0,r)$, we know that
\[ g_{i_{0}}\Big( G(z)+  C_{M-1}^{j}\cdot \Big(D(G)(z)\Big) (w)\Big)
= g(z)+C_{M-1}^{j} \cdot \Big( D(g)(z)\Big) (w) +
R\]
where $R=R_{i_{0}}\Big( C_{M-1}^{j}\cdot \Big( D(G)(z)\Big) (w), z \Big)$
and $j=1$ or $-1$, and 
\[ |R| \leq L_{0} C_{M-1}^{1+\alpha}||D(G)(z)||^{1+\alpha} |w|^{1+\alpha}.\]
But for $z_{0}=z$ and $z_{i}=g_{M-i}\circ \cdots \circ g_{i_{M}}(z)$,
$i=1$, $2$, $\ldots$, $M$, 
\[ ||D(G)(z)||= \prod_{1\leq k\leq M}||D(g_{i_{k}})(z_{M-k})|| \leq
\prod_{1\leq k\leq M}l_{i_{k}}(z_{M-k}).\]
Hence, by $K_{i}(z) \leq \Big( 1/ l_{i}(z)\Big)^{\beta}$ for all $i$, 
we have that
\[ ||D(G)(z)||^{1+\alpha}\leq 
\Big( \prod_{1\leq k\leq M}s_{i_{k}}(z_{M-k}) \Big) l^{(\alpha -\beta )M}.\]
Let $B_{M}= \Big( L_{0}/s \Big) C_{M-1}^{1+\alpha} 
\delta^{\alpha} l^{(\alpha -\beta )M}$, then
\[ |R|
\leq B_{M} \Big(
\prod_{0\leq k\leq M}s_{i_{k}}(z_{M-k}) \Big) |w|.\]
Since $B_{M} \leq \Theta_{\varepsilon} \delta^{\beta} l^{(\alpha
-\beta )M} \leq \delta^{\beta} l^{(\alpha -\beta )M}$, we get that 
$C_{M-1} + B_{M} \leq C_{M}$. Now we can conclude from the estimates that 
$g(w)-g(z)$ is in  
$C_{M}\cdot
\Big( D(g)(z)\Big) \Big( B(0, r)\Big)$ and if $|w|=r$,
$g(w)-g(z)$ is outside of $C_{M}^{-1}\cdot \Big( D(g)(z)\Big) \Big( B(0,
r)\Big)$. The proof of the claim is completed.

Take $C=C_{\infty}(\varepsilon)$. Then $\delta$ and $C$ 
are the functions we want. This completes the proof of Theorem 1.

The proof of Corollary 1 is similar.

\vskip10pt
\noindent {\bf \S 3 Proof of Theorem 2.}

According to Theorem 1, each $g_{w_{p}}(\overline{W})$  
contains a translation of the ellipse $C^{-1}\cdot \Big( D(g_{w_{p}})(z) \Big) \Big( B(0,
1)\Big)$ and 
is contained in a
translation of the ellipse $C\cdot \Big( D(g_{w_{p}})(z)\Big) \Big( B(0, 1)\Big)$ 
where $C$ 
is independent of $w_{p}$ and $z$.
For every $w_{p}=i_{0}i_{1}\cdots i_{p-1}$ in $\Sigma_{p}$, let
$g_{w_{p}}=g_{i_{0}}\circ g_{i_{1}}\circ \cdots \circ g_{i_{p-1}}$. 
Since all $g_{i}$ are contracting, there is a constant $0< \lambda_{0}<1$
such that $diam\Big( g_{w_{p}}(\overline{W})\Big) \leq \lambda_{0}^{p}$ for all
$w_{p}\in \Sigma_{p}$.  Thus 
$\{ g_{w_{p}}(\overline{W}); w_{p} \in \Sigma_{p}\}$ is a cover of
$\Lambda$ for every $p$ and 
$\tau_{p} = \max \{ diam(g_{w_{p}}(\overline{W}));  w_{p} \in
\Sigma_{p} \}$ tends to zero as $p$ tends to $\infty$. 
Use Theorem 1 again, the Hausdorff dimension \cite{f} of 
$\Lambda $ is a unique number $t_{0} >0$
satisfying 
\[ \lim_{p\mapsto \infty} \sum_{w_{p} \in \Sigma_{p} }\Big(
diam \Big( g_{w_{p}}(\overline{W})\Big) \Big)^{t} =\infty \hskip7pt for \hskip7pt 
t< t_{0} \hskip7pt and \] 
\[ \lim_{p\mapsto \infty} \sum_{w_{p} \in \Sigma_{p}  }\Big(
diam \Big( g_{w_{p}}(\overline{W})\Big) \Big)^{t} =0 \hskip7pt for \hskip7pt t> t_{0}.
\]

Let $l_{w_{p}}(z)$ and $s_{w_{p}}(z)$ be the lengths of longest and shortest
axes of 
the ellipse $\Big( D(g_{w_{p}})\Big) \Big( B(0,1)\Big)$. Then we have that
\[ C^{-1} \cdot s_{w_{p}}(z) \leq diam\Big( g_{w_{p}}(\overline{W})\Big) \leq C \cdot l_{w_{p}}(z).\]
One of the crucial points is that 
\[ l_{w_{p}}(z)\leq l_{i_{0}}(z_{p-1})\cdots l_{i_{p}}(z_{0}) \hskip7pt
and \hskip7pt s_{w_{p}}(z)\geq s_{i_{0}}(z_{p-1})\cdots s_{i_{p}}(z_{0})\]
where $z_{k}=g_{i_{p-k}}\circ \cdots \circ g_{i_{p-1}}(z)$.
Because of these two inequalities, we can conclude our proof 
by Gibbs theory (see, for example, \cite{b,r,s}) 
as follows: for any $t>0$, 
\[ \Big( diam\Big( g_{w_{p}}(W)\Big) \Big)^{t} \leq C_{1}\cdot \exp
\Big( \sum_{k=0}^{p-1} t \phi_{up} (w^{k})\Big)\hskip7pt and \]
\[ \Big( diam\Big( g_{w_{p}}(W)\Big) \Big)^{t} \geq C^{-1}_{1} \cdot \exp
\Big( \sum_{k=0}^{p-1} t \phi_{lo} (w^{k})\Big)\]
where $\pi (w^{k})=z_{k}$ and $C_{1}$ is a constant.
Suppose $\mu_{t_{up}\phi_{up}}$ and $\mu_{t_{lo}\phi_{lo}}$ are the Gibbs
measures of $t_{up}\phi_{up}$ and $t_{lo}\phi_{lo}$ on $\Big( \Sigma_{\infty},
\sigma \Big)$. 
Because $P(t_{up}\phi_{up})=0$ and $P(t_{lo}\phi_{lo})=0$, there is a constant
$d>0$ such that  
\[ \mu_{t_{up}\phi_{up}} (\Lambda_{w_{p}}) \in [d^{-1}, d] 
\exp \Big( \sum_{k=0}^{p-1}t_{up}\phi_{up} \Big( \sigma^{\circ
k}(w_{0})\Big) \Big) \hskip7pt
{\rm and} \]
\[ \mu_{t_{lo}\phi_{lo}} (\Lambda_{w_{p}}) \in [d^{-1}, d] 
\exp \Big( \sum_{k=0}^{p-1}t_{lo}\phi_{lo} \Big( \sigma^{\circ
k}(w_{0})\Big) \Big) \]
where $w_{0}\in \Lambda_{w_{p}}=\{ w\in \Sigma;  w=w_{p}\cdots \}$.
Hence there is a constant $C_{2}>0$ such that 
\[ \Big( diam\Big( g_{w_{p}}(W)\Big) \Big)^{t_{up}} \leq C_{2}\cdot
\mu_{t_{up}\phi_{up}}(\Lambda_{w_{p}}) \hskip7pt
{\rm and} \]
\[ \Big( diam\Big( g_{w_{p}}(W)\Big) \Big)^{t_{lo}} \geq C^{-1}_{2}\cdot
\mu_{t_{lo}\phi_{lo}}(\Lambda_{w_{p}}) . \]
Moreover,
\[ \sum_{w_{p} \in \Sigma_{p}} \Big( diam\Big( g_{w_{p}}(W)\Big)
\Big)^{t_{up}} \leq C_{2}\cdot \sum_{w_{p}\in \Sigma_{p} }
\mu_{t_{up}\phi_{up}}(\Lambda_{w_{p}})= C_{2}\hskip7pt and\]  
\[ \sum_{w_{p}\in \Sigma_{p} } \Big( diam\Big( g_{w_{p}}(W)\Big)
\Big)^{t_{lo}} \geq C^{-1}_{2} \cdot \sum_{w_{p}\in \Sigma_{p} }
\mu_{t_{lo}\phi_{lo}}(\Lambda_{w_{p}})= C^{-1}_{2}. \]
This implies that $t_{lo} \leq HD(\Lambda) \leq t_{up}$.
The proof is completed.

\vskip5pt
{\bf Proof of Corollary 2.} For $\phi=\phi_{lo,\lambda}$ (or 
$\phi_{up,\lambda}$),
the inverse of $P(t\phi)$ is continuous on $P$ and $\lambda$. This implies that
$t_{lo,\lambda}$ (or $t_{up, \lambda}$) tends to $t_{lo,\lambda_{0}}$ 
(or $t_{up,\lambda_{0}}$) as $\lambda$
goes to $\lambda_{0}$. But, $t_{lo,\lambda_{0}}=t_{up,\lambda_{0}}=
HD(\lambda_{0})$ because all
$g_{i,\lambda_{0}}$ are conformal. This completes the proof.

{\bf Proof of Corollary 3.} Let $\lambda =( b, c)$ (or $\lambda =(\gamma -n,
c)$) and $|\lambda| = |b|+|c|$ (or $|\lambda | =|\gamma-n|+|c|$). There
is a neighborhood $W$ of $S^{1}=\{ z\in {\bf C}; |z|=1\}$ so that 
$f$ is expanding on $\overline{W}$ for small $|\lambda|$. 
Let $g_{0, \lambda}$, $\ldots$, $g_{n-1, \lambda}$ be the 
inverse branches of $f|\overline{W}$. Then 
${\cal G}_{\lambda}$, the semigroup generated by 
$g_{0,\lambda}$, $\ldots$, $g_{n-1,\lambda}$, is regular and Markov for
$\lambda$ with small $|\lambda|$. Now the proof follows from Corollary 2
because for each
$\lambda = (0,c)$ with small $|c|\neq 0$, all $g_{i,
\lambda}$ are 
conformal and the Hausdorff dimension $HD(\lambda)$ of $J_{\lambda}$ is
greater than one.   

\vskip10pt
\noindent {\bf \S 4 Higher dimensional regular semigroups and some remarks.}

Suppose ${\bf E^{m}}$ is the $m$-dimensional Euclidean space, 
$V_{i}\subset U_{i}$, $i=0$, $\ldots$, $n-1$, are pairs of open
sets of ${\bf E}^{m}$ with $\overline{V_{i}}\subset U_{i}$ and $f_{i}$
from $\overline{V}_{i}$ onto $\overline{W}_{i}$ are
$C^{1+\alpha}$ diffeomorphisms such that the inverses $g_{i}$ of $f_{i}|\overline{V}_{i}$ are
contracting. 
Let ${\cal G}_{m}=\langle g_{0}, g_{1}$, $\ldots$, $g_{n-1} \rangle $ be 
the semigroup
generated by all $g_{i}$. Then $l$ and $K$ for ${\cal G}_{m}$ can be 
defined similarly. Again ${\cal G}_{m}$ is said to be regular if 
$K< 1/l^{\alpha}$.
Let $B(x,r)$ be the closed ball of radius $r$ centered at $x$ of ${\bf E}^{m}$.
The higher dimensional version of Theorem 1 is the following: 

{\bf Theorem 3} (geometric distortion). {\em Suppose
${\cal G}_{m}=\langle g_{0}$, $g_{1}$, $\ldots$, $g_{n-1}\rangle $ is regular. 
There are two functions
$\delta=\delta (\varepsilon)>0$ and $C=C(\varepsilon)\geq 1$ with 
$\delta (\varepsilon) \mapsto 0$ and $C(\varepsilon)\mapsto 1$ as $\varepsilon
\mapsto 0+$ such that 
\[ g\Big(
B(x, r)\Big)\supset g(x)+C^{-1}\cdot \Big( D(g)(x)\Big) \Big( B(0, r )\Big)
\hskip5pt and \]
\[ g\Big(
B(x, r)\Big) \subset
g(x)+C\cdot \Big( D(g)(x)\Big)\Big( B(0, r )\Big) \]
for any $0< r \leq \delta (\varepsilon)$, any $g\in
{\cal G}_{m}$ and any $x\in \overline{W}$.}

{\bf Remark 2.} Similarly, we have the higher dimensional versions of 
Corollary 1 and Theorem 2. 
We learned recently that Gu
\cite{gu} showed another upper bound (in higher dimensional case) 
which is similar to that in Theorem 2.

{\bf Remark 3.} Suppose $f_{\lambda}(z) =
z^{2}|z|^{(\gamma-2)} + c$ where $\lambda =(\gamma -2,
c)$. From Corollary 3 and Remark 1, there is an interesting
picture on the parameter space $\lambda$ (three dimensional space) 
near the point $(0, 0)$: there are small sectors $T_{1}$ and $T_{2}$ (see
Fig. 4) such
that for $\lambda $ in $T_{1}$, $J_{\lambda}$ is a smooth circle and for
$\lambda $ in $T_{2}$, $J_{\lambda}$ is a fractal circle with Hausdorff
dimension $>1$. From computer pictures of $J_{\lambda}$ for small $|\lambda
|$, we conjecture that there is a topological surface $S$ passing $(0,0)$
in a small ball centered at $(0,0)$
such that in the right hand side of $S$, $J_{\lambda}$ is a smooth circle
and in the left hand side of $S$ (but not on the $(\gamma-2)$-axis), 
$J_{\lambda}$ is a fractal circle with Hausdorff dimension $>1$ (see
Fig. 5). We may call $S$ the {\em boundary of fractalness}. If $S$ exists, 
what can be said about its shape ?

\vskip5pt
\centerline{\psfig{figure=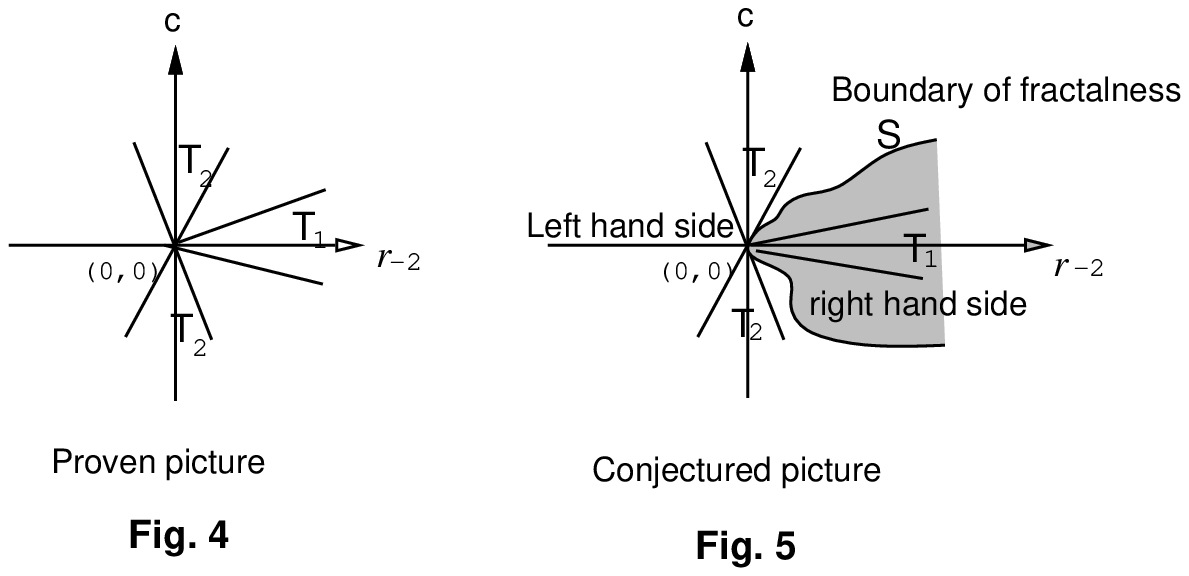}}

\vskip5pt
{\bf Remark 4.} Sullivan \cite{s} has considered 
quasiconformal deformations of
analytic and expanding systems and Gibbs measures. Moreover, he
also studied (uniform) quasiconformality in geodesic flows of 
negatively curved manifolds. One wonders if Theorem 1 can be used to extend
some results \cite{s} to non-conformal expanding systems (or hyperbolic
systems) with the compatibility condition
$K< 1/l^{\alpha}$ and to geodesic flows of 
negatively curved manifolds with pinched condition.

\end{document}